\documentclass[psamsfont]{amsart}

\usepackage{graphics}
\usepackage[dvips]{epsfig}
\usepackage{psfrag}
\usepackage{pinlabel}
\usepackage{amsmath}
\usepackage{amsfonts}
\usepackage{latexsym}
\usepackage{amssymb}
\usepackage[usenames]{color}

\input xy
\xyoption{all}

\newtheorem{defn}{Definition}[section]
\newtheorem{Thm}[defn]{Theorem}

\newtheorem{Rem}[defn]{Remark}

\newtheorem{Question}[defn]{Question}

\DeclareMathAlphabet{\mathdj}{U}{msb}{m}{n}
\title{Some non-collarable slices of Lagrangian surfaces.}
\author{Baptiste Chantraine}
\address{D\'epartement de Math\'ematiques\\
Universit\'e Libre de Bruxelles\\
CP 218\\
Boulevard du Triomphe\\
1050 Bruxelles\\   
 Belgique.}
\email{bchantra@ulb.ac.be}

\thanks{The author is supported by a Mandat Charg\'e de recherche from the Fonds National de la Recherche Scientifique, Communaut\'e fran\c caise de Belgique. The author also acknowledge the help of the research networking programme CAST funded by the European Science Foundation.}

\begin{document}
\maketitle
\begin{abstract}
  In this note we define the notion of collarable slices of Lagrangian submanifolds. These are slices of Lagrangian submanifolds which can be isotoped through Lagrangian submanifolds to a cylinder over a Legendrian embedding near a contact hypersurface. Such a notion arises naturally when studying intersections of Lagrangian submanifolds with contact hypersurfaces. We then give two explicit examples of Lagrangian disks in $\mathbb{C}^2$ transverse to $S^3$ whose slices are non-collarable.
\end{abstract}

\section{Introduction}
\label{sec:introduction}

In this note we raise the question of finding collarable slices of Lagrangian submanifolds in the symplectisation of a co-oriented contact manifold $(Y,\xi=\ker \alpha)$. We then provide examples of Lagrangian disks whose boundary are non-collarable. The symplectisation of $(Y,\xi)$ is the manifold $\mathbb{R}\times Y$ with symplectic form $d(e^t\alpha)$ where $t$ parametrizes the $\mathbb{R}$ direction. We begin by defining the object.
\begin{defn}
  Let $L\subset\mathbb{R}\times Y$ be a Lagrangian submanifold. We say that the slice $L\pitchfork \big(\{t_0\}\times Y\big)$ is collarable if there exists $\varepsilon_0>0$ and a Lagrangian isotopy $L_s$ of $L$ such that:
  \begin{itemize}
\item {$\forall s\in [0,1]\, L_s=L$ outside $[t_0-\varepsilon_0,t_0+\varepsilon_0]\times Y$}
\item {$\forall s\in[0,1]\,L_s\pitchfork\big(\{t_0\}\times Y\big)$}
\item {$\exists 0<\varepsilon<\varepsilon_0$ such that $\frac{\partial}{\partial t}$ is tangent to $L_1\cap\big( (t_0-\varepsilon,t_0+\varepsilon)\times Y\big)$.}
  \end{itemize}
\end{defn}

The second condition ensures that the smooth type of the slice does not change along the isotopy and the third asserts that locally $L_1$ is a cylinder over a Legendrian embedding $\Lambda$ in $(Y,\xi)$.

\begin{Rem}
  Note that, by the first condition, the definition is purely local and we do not ask for any constraint on $L$; it could be closed, non-compact, non-orientable etc. We only require the intersection with the contact slice to be transverse.
\end{Rem}

The notion of collarable Lagrangians arises in two natural contexts. The first is in the context of closed Lagrangian submanifolds $L$ in a symplectic manifold $(M,\omega)$. Let $Y$ be a (separating for simplicity) contact hypersurface of $M$ with transverse Liouville vector field $V$. The flow of the Liouville vector field allows us to identify a neighborhood of $Y$ with a part of the symplectisation of $(Y,\iota_V\omega)$. If, under this identification, the slice $L\pitchfork Y$ is collarable then cutting along $Y$ leads to two fillings of a Legendrian submanifold of $Y$. Such fillings were introduced and studied by the author in \cite{ChantraineConcordance}. From there one could deduce some topological constraints on both $L$ and the Legendrian slice. Stretching such a decomposition might also lead to some Mayer-Vietoris type decomposition of the Lagrangian Floer homology of $L$, however such considerations are far away from the purpose of this note.

The second one is the context of Lagrangian cobordisms, also introduced by the author in \cite{ChantraineConcordance}. We recall here the definition of this notion.
\begin{defn}
  Let $\Lambda^-$ and $\Lambda^+$ be two Legendrian submanifolds of $(Y,\xi)$. A Lagrangian cobordism from $\Lambda^-$ to $\Lambda^+$ is Lagrangian submanifold $L$ of $\mathbb{R}\times Y$ such that
$$\exists T>0\text{ s.t. } L\cap \big( (-\infty,-T)\times Y\big)=(-\infty,-T)\times \Lambda^-\text{ and }L\cap\big( (T,\infty)\times Y\big)=(T,\infty)\times \Lambda^+.$$
\end{defn}

We restrict ourselves in this note to the contact manifolds $(S^3,TS^3\cap iTS^3)$ in $\mathbb{C}^2$ and $(\mathbb{R}^3,\ker(dz-ydx))$. As the latter is the one point compactification of the former (even as contact manifold), we indistinctly denote by $(Y,\xi)$ any of those two. It is also convenient to recall that $(\mathbb{C}^2\setminus\{0\},\omega_0)$ is the symplectisation of the contact $S^3$ (with projection to $\mathbb{R}$ given by by $\ln (d(p)^2)$ where $d$ is the ``distance from the origin'' function), and that $T^*(\mathbb{R}_+\times\mathbb{R})$ is symplectomorphic to the symplectisation of $\mathbb{R}^3$ via the map $(q_1,p_1,q_2,p_2)\rightarrow (q_2,q_1p_2,ln(q_1),p_1)$. We denote by $S^3_\varepsilon$ (resp. $D^4_\varepsilon$) the sphere (reps. disk) of radius $\varepsilon$ in $\mathbb{C}^2$. 

\begin{figure}[ht!]
\labellist
\small\hair 2pt
\pinlabel {$1$} [bl] at 176 620
\pinlabel {$1'$} [bl] at 401 620
\pinlabel {$2$} [bl] at 118 489
\pinlabel {$2'$} [bl] at 451 489
\pinlabel {$3$} [bl] at 298 344
\pinlabel {$4$} [bl] at 298 179
\pinlabel {$5$} [bl] at 298 38
\pinlabel {$\emptyset$} [bl] at 96 22
\endlabellist
  \centering
  \includegraphics[height=10cm]{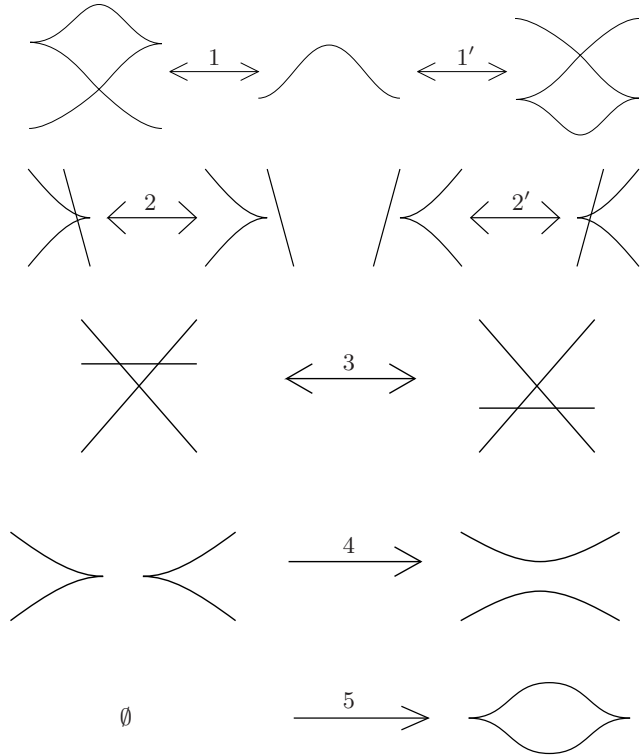}
  \caption{Front projection of simple moves along Lagrangian cobordisms.}
  \label{fig:Local_move}
\end{figure}

All constructions of explicit examples of Lagrangian cobordisms between Legendrian knots, to this day, involve sequences of simple moves of the front projections of Legendrian knots summarized in Figure \ref{fig:Local_move} (see \cite{ChantraineConcordance}, \cite{KalHoEk} and \cite{TraySabTQFT} ). One recovers the original simple move by cutting the cobordism along some collection of collarable slices, we give a name to Lagrangian cobordisms arising this way.
\begin{defn}
  A Lagrangian cobordism $L$ from $\Lambda^-$ to $\Lambda^+$ is decomposable if there exists $s_1<s_2<\cdots<s_k\in\mathbb{R}$ such that:
  \begin{itemize}
  \item {Each slice $L\cap(\{s_i\}\times Y)$ is collarable}
\item {$L\cap\big( (-\infty,-s_1)\times Y\big)=(-\infty,-s_1)\times\Lambda^-$ and $L\cap\big( (s_k,\infty)\times Y\big)=(s_k,\infty)\times\Lambda^+$}
\item {The fronts of the slices $L\cap\{s_i\}\times Y$ and $L\cap\{s_{i+1}\}\times Y$ differ only by one of the events shown in Figure \ref{fig:Local_move} and $L\cap[s_i,s_{i+1}]\times Y$ is the elementary Lagrangian cobordism realising the bifurcation.}
  \end{itemize}
\end{defn}

The fact that only such decomposable cobordisms were constructed raises the following question:

\begin{Question}\label{sec:que1}
  Is any exact Lagrangian cobordism between two Legendrian links in $(Y,\xi)$ decomposable?
\end{Question}
 A sub-question of the preceding, whose positive answer would be the first step toward a positive answer to Question \ref{sec:que1}, would be:

 \begin{Question}\label{sec:que2}
   Is any transverse slice of an exact Lagrangian cobordism collarable?
 \end{Question}

A second step toward a positive answer to Question \ref{sec:que1} would be to classify the possible critical points of the projection of a Lagrangian cobordism to $\mathbb{R}$.

In the present note we exhibit some non-collarable Lagrangian slices arising as boundary of (exact) Lagrangian disks in $\mathbb{C}^2$. Namely we prove:

\begin{Thm}\label{thm:1}
  There exists a Lagrangian cobordism $T$ from the trivial Legendrian knot $K_0$ to the empty set with one global maximum $m$ for the ``distance from the origin'' function such that for any sufficiently small $\delta>0$ the slice $T\cap S^3_{m-\delta}$ is non-collarable.
\end{Thm}

The part of $T$ which lies on the concave side of $S^3_{m-\delta}$ is thus a Lagrangian disk (hence exact) which provides a local model for a maximum of the projection to $\mathbb{R}$ restricted to a Lagrangian submanifold in the symplectisation. We also provide an example of a Lagrangian disk in the convex part of the symplectisation, namely:
\begin{Thm}\label{thm:2}
  There exist Lagrangian disks in $D^4$ transverse to $S^3$ whose intersections with $S^3$ are non-collarable.
\end{Thm}

It is interesting to compare the existence of these disks with the result of \cite{EliPollag} where it is shown that any Lagrangian disk whose boundary is the trivial Legendrian knot $K_0$ is equivalent to the flat Lagrangian disk.

So far none of these examples provide negative answer to Question \ref{sec:que2}. Indeed none of these (exact) Lagrangian disks appear here as part of an exact Lagrangian cobordism. The example of Theorem \ref{thm:1} emphasizes the importance of the exact condition in Question \ref{sec:que2} as the Lagrangian cobordism $T$ is explicitly non exact. If the answer to Question \ref{sec:que2} were positive such Lagrangian disks could not be part of an exact Lagrangian cobordism.

In Section \ref{sec:exist-obstr-lagr}, we recall the previous results on the construction of Lagrangian cobordisms of Figure \ref{fig:Local_move}. We also recall the main obstruction to the existence of Lagrangian fillings from \cite{ChantraineConcordance}. In Section \ref{sec:local-maximum} and  \ref{sec:knott-lagr-disk}, we give the constructions of the two Lagrangians of Theorems \ref{thm:1} and \ref{thm:2} and prove that the slices are non-collarable.

\section*{acknowledgements}\label{ackref}
It was during the focused week on generating families of the special semester on contact and symplectic topology in Nantes, in 2011, that I had the opportunity to talk about these constructions where they appear to be of some interest. I would like to thank the organizers of the semester Vincent Colin, Paolo Ghiggini and Yann Rollin. I also thank Margherita Sandon for organizing this focused week. Finally I thank Petya Pushkar for suggesting me to write this note.

\section{Existence and obstruction to Lagrangian cobordism.}
\label{sec:exist-obstr-lagr}

We begin by describing the construction of the elementary cobordisms realizing the moves of Figure \ref{fig:Local_move} motivating Question \ref{sec:que1}.

The first three moves are Legendrian Reidemeister moves arising along generic Legendrian isotopies. The existence of Lagrangian cobordisms realizing those moves follows from \cite[Theorem 1.2]{ChantraineConcordance}.

Move number $5$ arises from the standard capping of the trivial Legendrian knot $K_0$ in $S^3$ or $\mathbb{R}^3$. This capping in  $S^3$ is simply the intersection of the Lagrangian plane $\mathbb{R}^2\subset\mathbb{C}^2$ with the $4$-ball $D^4$. The capping in $\mathbb{R}^3$ can be deduced via the contactomorphism from $S^3/N$ to $\mathbb{R}^3$. Using the symplectomorphism $\mathbb{R}\times\mathbb{R}^3\simeq T^*(\mathbb{R}_+\times\mathbb{R})$ it can also be described by the generating family $F(s,t,\eta)=s\big(\frac{\eta^3}{3}-\frac{3}{2}(\rho(s)-t^2)\eta\big)$, where $\rho$ is a smooth increasing function equal to $-1$ on $(0,\delta]$ and $1$ for $s>1-\delta$.

\begin{figure}[ht!]
  \centering
  \includegraphics[height=4cm]{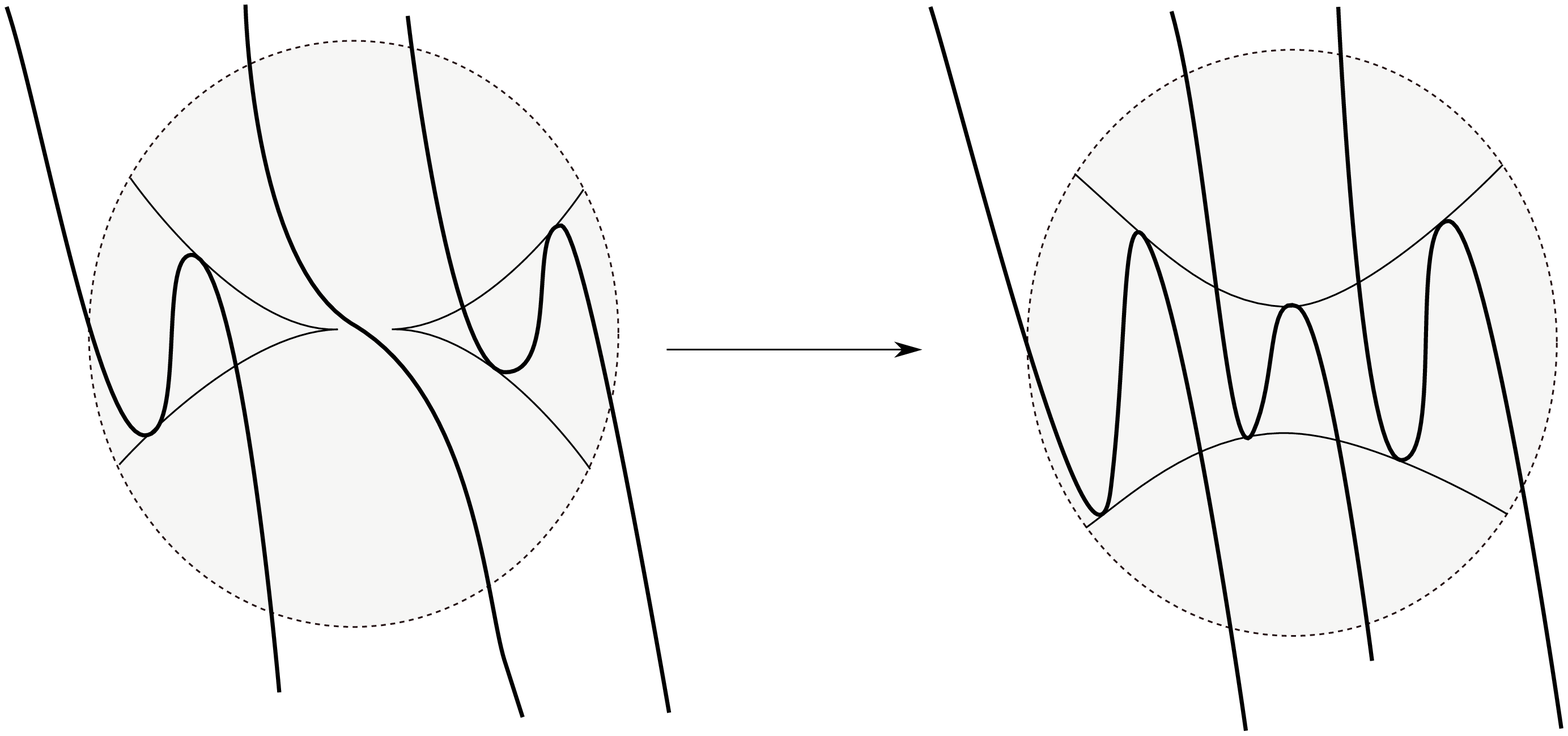}
  \caption{Generating family for an index $1$ critical point of the projection to $\mathbb{R}$.}
  \label{fig:GenFam}
\end{figure}

Move number $4$ was announced by T. Ekholm, K. Honda and T. K\'alm\'an in \cite{KalHoEk}. J. Sabloff and L. Traynor in \cite{TraySabTQFT} realize this move in the context of Lagrangian cobordisms with compatible generating families. The move can be realized by interpolating between the two generating families of Figure \ref{fig:GenFam}.

We recall now the main obstruction we use here to the existence of a Lagrangian filling of a Legendrian knot.
\begin{Thm}[\cite{ChantraineConcordance}]\label{thm:filling}
  Let $\Sigma$ be an oriented Lagrangian submanifold of $D^4$ such that $\partial \Sigma$ is a Legendrian submanifold $K$ of $S^3$. Then the following holds:
$$r(K)=0,$$
$$tb(K)=2g(\Sigma)-1=TB(\mathcal{L}(K)),$$
$$g(\Sigma)=g_s(\mathcal{L}(K)),$$
where $\mathcal{L}(K)$ is the smooth isotopy type of $K$, $TB(\mathcal{L}(K))$ is the maximal Thurston-Bennequin number among Legendrian representatives of $\mathcal{L}(K)$ and $g_s(\mathcal{L}(K))$ is the $4$-ball genus of $\mathcal{L}(K)$.
\end{Thm}

\section{The local maximum.}
\label{sec:local-maximum}

\begin{figure}[ht!]
\labellist
\small\hair 2pt
\pinlabel {$\gamma$} [bl] at 219 300
\pinlabel {$p$} [bl] at 341 395
\pinlabel {$q$} [bl] at 466 327
\pinlabel {$D$} [bl] at 300 300
\endlabellist
  \centering
  \includegraphics[height=6cm]{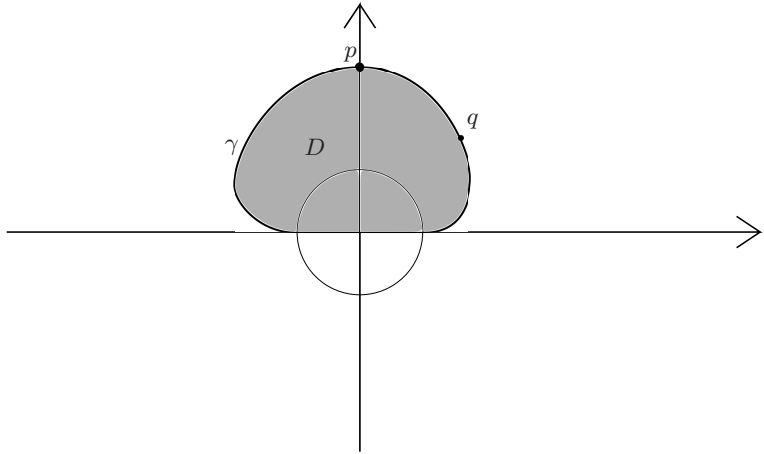}
  \caption{A simple close curve in $\mathbb{C}$.}
  \label{fig:simplecurve}
\end{figure}

\begin{figure}[ht!]
\labellist
\small\hair 2pt
\pinlabel {$S^3$} [bl] at 170 109
\pinlabel {$K_0$} [bl] at 203 308
\pinlabel {$T$} [bl] at 262 377
\pinlabel {$S^3_{m-\delta}$} [bl] at 100 454
\endlabellist
  \centering
  \includegraphics[height=4cm]{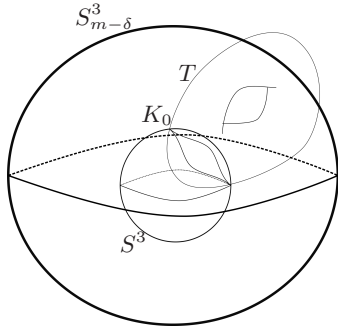}
  \caption{A Clifford torus in $\mathbb{C}^2$.}
  \label{fig:clifford}
\end{figure}

The disk of Theorem \ref{thm:1} arises from the Clifford torus in $\mathbb{C}^2$. Consider the simple closed curve $\gamma$ in $\mathbb{C}$ shown in Figure \ref{fig:simplecurve}. One only asks that the intersection of $\gamma$ with the unit disk is the interval $[-1,1]\times\{0\}$ in $\mathbb{R}^2\simeq\mathbb{C}$, and that the distance function restricted to $\gamma$ has one global maximum $\frac{m}{\sqrt{2}}$ at $p$. Then $\gamma\times\gamma\subset\mathbb{C}^2$ is a Lagrangian torus $T_0$ whose intersection with $D^4$ is the Lagrangian filling of the trivial Legendrian knot $K_0$. Now consider the slice $S^3_{m-\delta}\cap T$ (see Figure \ref{fig:clifford}). If the slice were collarable one would get a Legendrian knot in $S^3$ with a Lagrangian filling of genus $1$. From Theorem \ref{thm:filling}, the $4$-ball genus of this knot would be $1$, however it is clearly the boundary of a disk on the concave side, hence its $4$-ball genus is $0$.

The part of $T_0$ which lies outside of $D^4$ is a Lagrangian cobordism from $K_0$ to the empty set. It is not exact as $\{q\}\times D$ is a symplectic disk with boundary on $T$. As $\delta$ can be chosen arbitrarily small, the disk provides a local model for the maximum of the projection of a Lagrangian submanifold of the symplectisation to $\mathbb{R}$. 

\section{A knotted Lagrangian disk in $D^4$}
\label{sec:knott-lagr-disk}

The non-collarable Lagrangian slices of Theorem \ref{thm:2} arise from knotted Lagrangian disks in $D^4$. We construct them using work of L. Rudolph \cite{RudolphAlg}, where it is shown that every knot which arises as the closure of a quasi-positive braid is the intersection of a complex curve with $S^3$. Such knots are called \textit{quasi-positive}. We recall the definition of a quasi-positive braid.

\begin{defn}
  Let $B_n$ be the braid group with $\sigma_1\cdots\sigma_{n-1}$ its standard generators. An element of $B_n$ is quasi-positive if it can be written as
$$\sigma=w_1\sigma_{i_1} w_1^{-1}\cdots w_k\sigma_{i_k} w_k^{-1},$$
where the $w_i$'s are any element of $B_n$. 
\end{defn}

The piece of algebraic curve inside the $4$-ball, whose boundary is a quasi-positive knot, is described as a $n$ fold branched cover over the disk with $k$ simple branched points. The Riemann-Hurwitz formula gives that its Euler characteristic is $n-k$. It is thus a disk if $k=n-1$.

\begin{figure}[ht!]
  \centering
  \includegraphics[height=3cm]{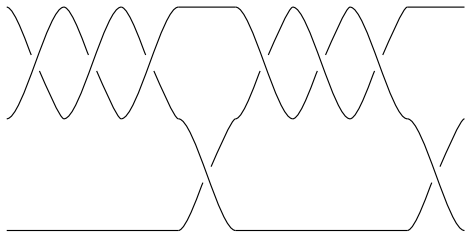}
  \caption{Quasi-positive braid whose closure yields the mirror of $8_{20}$.}
  \label{fig:braid}
\end{figure}

Note that any holomorphic disk in $D^4$ leads to a Lagrangian disk by applying the linear transformation $(x_1,y_1,x_2,y_2)\rightarrow(x_1,x_2,y_1,y_2)$. Now, in Rolfsen's knot table \cite{Rol}, the mirror of the knot $8_{20}$ is quasi-positive as it is the closure of $(\sigma_1^{-3})\sigma_2(\sigma_1^3)\sigma_2$ in $B_3$ (see Figure \ref{fig:braid}). In this case $n=3$ and $k=2$,  thus it arises as the boundary of a Lagrangian disk $D$ in $D^4$. However a computation due to L. Ng in \cite{NgTB} shows that the maximal Thurston-Bennequin of $K$ is $-2$. This implies that the slice of Lagrangian given by $S^3\cap D^4$ is non-collarable.

Note that quasi-positive knots of this type are numerous, we only mentioned the first one appearing in Rolfsen's knot table. Many of them have maximal Thurston-Bennequin number different from $-1$ and provide thus non-collarable slices. Also note that this example is of a different nature than the previous one. It cannot appear as a local minimum of the projection to $\mathbb{R}$. In fact, any slice of a Lagrangian submanifold near a local minimum of the projection is collarable. This can be seen viewing the Lagrangian near the minimum as the graph of an exact form $df$ on the (Lagrangian) tangent plane at the minimum. Perturbing $f$ to be constant near $0$ gives a deformation to a collar over the Legendrian trivial knot.

\bibliographystyle{plain}
 \bibliography{biblio}

\end{document}